\documentclass{ifacconf}

\usepackage{graphicx}      				% include this line if your document contains figures
\usepackage{natbib}        				% required for bibliography
\usepackage{amsmath,amssymb,amsfonts}	% align environment
\usepackage{siunitx} % murp added for SI units
\usepackage{color}

\newcommand{\U}{\mathcal{U}}
\newcommand{\M}{\mathcal{M}}
\newcommand{\DMM}{\left[\partial \mathcal{M}\right]_-}
\newcommand{\DMO}{\left[\partial \mathcal{M}\right]_0}
\newcommand{\Mcomp}{\mathcal{M}^{\textsf{C}}}
\newcommand{\A}{\mathcal{A}}

\newcommand{\DAM}{\left[\partial \mathcal{A}\right]_-}

\newcommand{\RR}{\mathbb{R}}

\newtheorem{definition}{Definition}

\begin{document}
\begin{frontmatter}

%\title{Sustainable Management of interconnected Food Production Systems via Theory of Barriers\thanksref{footnoteinfo}} 
% Application of Theory of Barriers for sustainable management of interconnected Food Production Systems
\title{Sustainability Analysis of Interconnected Food Production Systems via Theory of Barriers}
%\title{Robustness and Viability Analysis of Ecosystems for Sustainable Food Production via the Theory of Barriers\thanksref{footnoteinfo}} 
%1) on barrier theory for viable populations in ecosystems with interacting species.
%2) viable analysis of ecosystems for sustainable food production via theory of barriers
%3) viable population analysis for predator-prey systems/for food production via barrier theory
%4) viable analysis on ecosystems via barrier theory : A case study on controlled agriculture

\thanks[footnoteinfo]{This project has received funding from the European Social Fund (ESF).}

\author[First]{Tim Aschenbruck} 
\author[First]{Willem Esterhuizen} 
\author[First]{Murali Padmanabha}
%\author[Second]{Jean L\'{e}vine}
\author[First]{Stefan Streif}

\address[First]{Automatic Control and System Dynamics Laboratory, Technische Universit\"{a}t Chemnitz, 09111 Chemnitz, Germany (e-mail: tim.aschenbruck@etit.tu-chemnitz.de; willem.esterhuizen@etit.tu-chemnitz.de; murali.padmanabha@etit.tu-chemnitz.de stefan.streif@etit.tu-chemnitz.de).}
%\address[Second]{CAS, Math\'{e}matiques et Syst\`{e}mes, Mines-ParisTech, 35, rue Saint-Honor\'{e}, 77300 Fontainebleau,
%	France (e-mail: jean.levine@mines-paristech.fr).}
%\address[Third]{Electrical Engineering Department, 
%   Seoul National University, Seoul, Korea, (e-mail: author@snu.ac.kr)}

\begin{abstract}% Abstract of not more than 250 words.
%The interest for sustainable ecosystems has been increased in recent decades due to irreversible consequences based on humans. %18 words
%Overexploitation and extinction of species are the consequences of this actions. %11 words
%Strategies to prevent these dangers are viability approaches for a sustainable management of different ecosystems.
Controlled environment agriculture (CEA) is used for efficient food production. 
Efficiency can be increased further by interconnecting different CEA systems (e.g. plants and insect larvae or fish and larvae), using products and by-products of one system in the other.

%The overall system based on these interconnected systems can then be described by models of interacting species.

These interconnected systems define an overall system that can be described by models of interacting species.
	It is necessary to identify system parameters (e.g. initial species concentration, harvest rate, feed quality, etc.) such that the resources are not exhausted.
	For such systems with interacting species, modelled by the Lotka-Volterra equations, a set-based approach based on the recent results of the theory of barriers to exactly determine the so-called \emph{admissible set} (also known as viability kernel) and the \emph{maximal robust positively invariant set} is presented.	
	Using an example of a larvae-fish based production system, steps to obtain special trajectories which are the boundaries of the admissible set are shown.
	This admissible set is used to prevent the under and over population of the species in the CEA.
	Furthermore, conditions of the system parameters are stated, such that the existence of these trajectories can be guaranteed.

%We demonstrate our results based on an example from the literature to emphasize the simplicity of our analysis and an example based on a controlled environment for food production.
\end{abstract}

\begin{keyword}
%Five to ten keywords, preferably chosen from the IFAC keyword list.
Nonlinear systems, Continuous time systems, System analysis, Constraint satisfaction problems, Sets, Controlled environment agriculture 
\end{keyword}

\end{frontmatter}
%===============================================================================
\graphicspath{{./gfx/}}
\section{Introduction}
%{\color{red} 
%	\begin{itemize}
%		\item why is the sustainable foodproduction/ecosystem important for the future?
%		\item why do we use the predator prey model for this purpose?
%		\item why do we use set based methods? (a lot of uncertainties (measurements, parameters, etc.))
%		\item what have other people done? (viability analysis bayen, de lara, etc.)
%		\item outline
%	\end{itemize}
%}	

In recent years the general interest in sustainable ecosystems has increased, because of environmental changes and publicly well known challenges concerning, for example conservation of biodiversity, overexploitation of marine ecosystems, and exhaustion of fossil resources.
%Additionally, human actions influence the environment with irreversible consequences. For instance, the loss of species or natural resources.
It has also been recognised that human interventions yield ecological stresses with irreversible consequences. For instance, the loss of species or destruction of habitat.
Therefore, sustainable management of ecological systems becomes more and more important as described in \cite{deLara2008sustainableBOOK}.

Food production industry is one of the biggest contributors to this ecosystem imbalance \citep{Barnosky2011,Funabashi2018}.
	Modern farming techniques such as controlled environment agriculture (CEA), also used in space research, are evaluated to improve the productivity and address the sustainability aspects through the use of different organisms (e.g. plant, fishes, insects) interacting through resource exchange \citep{CONRAD2017}.
	Such studies are new in CEA and can utilize methods and tools developed for conventional farming and ecosystem management.

%%%%
%The dynamics of the organisms interacting either directly or through common resources---modelled predominantly using Lotka-Volterra equations---are widely used in the marine- and terrestrial ecosystem studies. 
%For example in sustainable management schemes for fisheries \citep{eisenack2006viabilityAnalysis} and sustainable production-consumption systems using an exhaustible natural resource \citep{martinet2007sustainabilityViableControl}, where the viability theory is applied to find the viability kernel to prevent exhaustion of species or resources.

The Lotka-Volterra equations for describing the species interaction and the viability theory for finding the viability kernel to prevent species exhaustion, are widely used in the marine- and terrestrial ecosystem studies.
For example, in \cite{bayen2019minimalTimeCrisisToViabKernel} the viability kernel was found to prevent the extinction of prey,
in \cite{eisenack2006viabilityAnalysis} the viability theory was applied to design management framework for fisheries,
and for exhaustible natural resource a viable control approach is presented via viability theory in \citep{martinet2007sustainabilityViableControl}.

 The aforementioned species interaction dynamics are also used to address different aspects of agriculture and food production. 
	Development of protocols for food production optimization \citep{FORT201723}, effects of organic management practices in pest control \citep{Schmidt2014}, 
	optimization of harvesting policy \citep{ZHANG2000}, design of ecosystems for nitrification in aquaponics system \citep{Graham2007}, and modeling the growth of \textit{Hermetia illucens} \citep{Anton2019} are some of the examples applying these interaction dynamics.

%% OLD
%Applying the viability theory and the widely studied Lotka-Volterra model on a food production system with species interaction to obtain viability kernel for sustainable production of interacting species is the main goal of this work.
%
%
%In this paper we analyse the \emph{admissible} set (also known as the viability kernel) and the \emph{maximal robust positively invariant} set, of a system described by the Lotka-Volterra equations with state and input constraints. 
%
%We use the theory of \emph{barriers in constrained nonlinear systems} as developed in \cite{DeDona_Levine2013barriers} and \cite{Ester_Asch_Streif_2019} to analyse these sets, describing how to construct them, and stating algebraic conditions of the system parameters under which they exist. 
%The usage of this theory simplifies the analysis drastically, which we show for a fish and insect-larvae example derived from a combination of aquaculture and insect farming under controlled environments.	
%%

In this paper we analyse the \emph{admissible set} and the \emph{maximal robust positively invariant} (MRPI) set of the widely studied Lotka-Volterra model, with state and input constraints, on a food production system with species interaction.
We use the theory of \emph{barriers in constrained nonlinear systems} as developed in \cite{DeDona_Levine2013barriers} and \cite{Ester_Asch_Streif_2019} to analyse these sets, describing how to construct them, and stating algebraic conditions of the system parameters under which they exist.
We use these sets to obtain information pertaining to management strategies for sustainable operation of food production systems.
We illustrate our results for a fish and insect-larvae example derived from a combination of aquaculture and insect farming under controlled environments.

%OUTLINE
The outline of this paper is as follows. In Section~\ref{sec:system} a brief overview of the model and its relationship with controlled environment agriculture is established. In Section~\ref{sec:2} we briefly cover the theory of barriers. 
Sections~\ref{sec:3} contains the main result of this paper with the analysis of the admissible and positively invariant set for the Lotka-Volterra equations. 
Section~\ref{sec_examples} presents a numerical examples and in Section~\ref{Sec_discussion} we discuss our analysis.
We conclude the paper in Section~\ref{sec_conclusion} with a summary and ideas for future research.

\section{Idea of sustainable interconnected controlled environment agriculture}\label{sec:system}
%\begin{itemize}
%	\item WHY FOR THIS APPLICATION? (CONTROLLED AGRICULTURE/FOOD PRODUCTION)
%	\item WHY PREDATOR-PREY?
%	\item WHY VIABILITY KERNEL?
%\end{itemize}

%In this section we present a general model for as system with interacting species known as predator-prey system. 
%This System is modelled by the Lotka-Volterra equations.
%Furthermore we explain the usage of this model in the field of sustainable controlled agriculture and in particular food production.
%For this purpose we use are interested in finding the admissible set (also known as viability kernel) for the model.

%In a CEA system producing single biomass (food) type, interaction can be established between the organism (species) growing and the resource supplied for its growth.
%	An example of such a system is the insect farming where interaction between the insect larvae and the substrate that it feeds on can be established by supplying insect feed to the system.
	In an interconnected CEA system producing different biomass (plant, fish, insect larvae) interaction between the connected systems is established through the exchange of products and by-products (e.g. plant waste as feed for larvae, larvae as feed for fish, fish waste as nutrients for plant).
	In case of a single CEA system the interconnection can be represented by the species interacting with a food source (e.g. larvae feeding on larvae feed, plant consuming nutrients etc.)

%	In a CEA system producing one biomass type, coupled with a system producing a second biomass type, interaction between the two systems is established through a common resource exchange.
%	For example, in aquaponics, interaction between a fish farm and a plant farm is established by cycling water between the two systems.
	
	Dynamics of such interactions between the species or between species and resources can be modelled using the Lotka-Volterra equations given by:
%\subsection{Model of Predator-Prey System}
%Consider the following Lotka-Volterra equations:
\begin{align}\label{eq_LotkaVolterra}
\begin{split}
\dot{x}_1 &= \alpha x_1 - \beta x_1 x_2 + u_1x_1, \\
\dot{x}_2 &= \delta x_1 x_2 - \gamma x_2 +  u_2x_2,
\end{split}
\end{align}
where $x_1\in\mathbb{R}_{\geq 0}$ denotes the number of prey (e.g. plant, larvae) and $x_2\in\mathbb{R}_{\geq 0}$ denotes the number of predators (e.g. larvae, fish).
The parameters $\alpha, \beta, \gamma$ and $\delta$ are positive constants. 
%and explained in the application framework of a food production example.
As per \cite{getz2012biomass}, the best applicable interpretation of the parameters of the Lotka-Volterra equations for biomass production are: $\alpha$, the intrinsic growth rate of the resource (prey);
$\beta$, the extraction rate per unit resource per unit consumer (predator); $\delta = \eta\beta$, where $\eta$ is the biomass conversion parameter; and $\gamma$, the intrinsic rate of decline of the consumer in the absence of resource.
%The positive constant parameters $\alpha$ and $\beta$ are the prey reproduction and eating rate; and $\gamma$ and $\delta$ are the predator mortality and reproduction-per-prey rate, respectively, as in \cite{royama1971PredPreyModels}. 
Furthermore, $u_i(t)\in[u_i^{\min}, u_i^{\max}]\subset\mathbb{R}$ for all $t$, with $u_i^{\min} \leq u_i^{\max}$ $i=1,2$, are constrained inputs that model the intervention into the predator-prey system.
In our framework of a CEA system the input is either the introduction of prey or the removal (harvest) of predators.
Additionally we have to introduce the limits for the number of prey and predators.
The lower limit is naturally given by zero, whereas for the goal of preventing extinction of a species, this number has to be positive.
We also introduce an upper limit, since, we act in a controlled environment with a given limit of the habitat for the species, thus, we also want to prevent overpopulation (overproduction).
Therefore we need to consider the following constraints
$x_1(t) \in [\underline{x}_1, \overline{x}_1] \; \forall t$, $ 0 \leq \underline{x}_1 < \overline{x}_1$,
$x_2(t) \in [\underline{x}_2, \overline{x}_2] \; \forall t$, $ 0 \leq \underline{x}_2 < \overline{x}_2$.
These constraints describe the limits for the under- and overpopulation.

%\subsection{Sustainability based on viability}
As mentioned in the introduction, viability theory plays a decisive role for sustainable strategies related to different ecosystems.
In particular, the identification of the viability kernel for predator-prey systems yields information to deduce harvesting strategies which prevent the extinction of a species.
In detail, this set describes the population number of predators and prey such that there exists at least one intervention strategy that results in the population remaining within the defined limits.
On the boundary of the admissible set there exist \emph{exactly one} intervention strategy that prevents under- or over-population.
The MRPI is interpretable as the population of predators and prey for which every intervention is sustainable.
In other word, there does not exist the possibility of crossing the limits for overpopulation or underpopulation.
In contrast to the MRPI, every population number \emph{outside} of the admissible set will lead to a violation of the constraints and therefore to an inevitable overpopulation or underpopulation.

\section{Summary of the Theory of Barriers} \label{sec:2}
This section presents a succinct review of the theory of barriers as presented in \cite{DeDona_Levine2013barriers} and \cite{Ester_Asch_Streif_2019}.

The considered nonlinear system subjected to state and input constraints is defined as follows:
\begin{align}
\dot{x}(t) & = f(x(t),u(t)),\quad x(t_0) = x_0,\quad u  \in \mathcal{U},\label{sys_eq_1}\\
g_i(x(t)) &\leq 0,\forall t\in[t_0,\infty[,\,\,i=1,2,\dots,p,\label{sys_eq_2}
\end{align}
where $x(t)\in\mathbb{R}^n$ denotes the state; $x_0$ is the initial condition at the initial time, $t_0$; $u(t)\in\mathbb{R}^m$ denotes the input, and $\mathcal{U}$ is the set of all Lebesque measurable functions that map $[t_0,\infty[\subset\mathbb{R}$ into $U\subset \mathbb{R}^m$, with $U$ compact and convex. The functions $g_i$, $i=1,\dots,p$, are constraint functions imposed on the state. We impose the same assumptions on the problem data as those specified in \cite{DeDona_Levine2013barriers} and \cite{Ester_Asch_Streif_2019}. Briefly, we assume that the functions $f$ and $g_i$ for $i=1,\dots,p$ are $C^2$ with respect to their arguments on appropriate open sets; that all solutions of the system remain bounded on finite intervals; and that the set $\{f(x,u):u\in U\}$ is convex. 
To lighten our notation, we introduce the following definitions:
\begin{align}
G\triangleq& \{x : g_i(x)\leq 0, \quad  \forall i\in\{1,2,...,p\}\},\label{def_G}\\
G_-\triangleq& \{x : g_i(x) <0, \quad  \forall i\in\{1,2,...,p\}\},\nonumber\\
G_0\triangleq& \{x : \exists i\in\{1,2,...,p\} \;\; \text{s.t.} \;\; g_i(x) =0\}.\nonumber
\end{align}
The set of indices of active constraints at $x$ is denoted by $\mathbb{I}(x)=\{i: g_i(x) = 0\}$.
Furthermore, we denote by $L_fg(x,u)$ the Lie derivative of a continuously differentiable function $g:\mathbb{R}^n \rightarrow \mathbb{R}$ with respect to $f(x,u)$ at the point $x$. Thus, $L_fg(x,u)=Dg(x)f(x,u)$.
The boundary of a set $S$ is denoted by $\partial S$. %, and its complement by $S^{\mathsf{C}}$. 
By $\mathbb{R}_{\geq 0}$ we refer to the set of nonnegative real numbers. In the following definitions, let $x^{(u,x_0,t_0)}(t)$ denote the solution to the differential equation at time $t$, initiating at $x_0$ at time $t_0$, with input $u\in \U$.

\begin{definition}
	The \emph{admissible set} of the system \eqref{sys_eq_1}-\eqref{sys_eq_2}, denoted by $\A$, is the set of initial states for which there exists a $u\in \U$ such that the corresponding solution satisfies the constraints \eqref{sys_eq_2} for all future time.
	\begin{align*}
	\mathcal{A} \triangleq \left\{ x_0 \in \RR^n: \exists u\in\mathcal{U}, \;\;  x^{(u,x_0,t_0)}(t)\in G \;\; \forall t \in [t_0,\infty[ \;\right\}.
	\end{align*}
\end{definition}
\begin{definition}
	A set $\Omega\subset\RR^n$ is said to be a \emph{robust positively invariant set} (RPI) of the system \eqref{sys_eq_1} provided that $x^{(u,x_0,t_0)}(t)\in\Omega$ for all $t\in[t_0,\infty[$, for all $x_0\in\Omega$ and for all $u\in\U$.
\end{definition}
\begin{definition}
	The \emph{maximal robust positively invariant set} (MRPI) of the system \eqref{sys_eq_1}-\eqref{sys_eq_2} contained in $G$, is the union of all RPIs that are subsets of $G$. Equivalently\footnote{shown in \cite{Ester_Asch_Streif_2019}, Proposition~2}
	%\cite[Proposition~2]{Ester_Asch_Streif_2019}},
	\begin{align*}
	\mathcal{M} \triangleq \left\{ x_0\in\RR^n: x^{(u,x_0,t_0)}(t)\in G, \,\,\forall u\in\mathcal{U}, \;\;  \forall t \in [t_0,\infty[ \;\right\}.
	\end{align*}
\end{definition}
%In general we distinguish between three different behaviours of the dynamical system depending on the state location in $\M$, $\A$, and $\Acomp$.
%Briefly, all trajectories initiating from $\M$ will stay within this set for all controls and will never violate the constraints, whereas all trajectories initiating from $\Acomp$ will definitely violate a constraint in the future, regardless of the control.
%The set $\A$ describes the situation where a future constraint violation depends on the chosen control function.

We introduce:
$\DAM \triangleq  \partial\mathcal{A}\cap G_-$ called the \emph{barrier} and $\DMM \triangleq \partial\mathcal{M}\cap G_-$,
called the \emph{invariance barrier}. These parts of the sets' boundaries are special, in that for any initial condition located on the barrier (resp. invariance barrier) there exists an input such that the resulting curve remains on the barrier (resp. invariance barrier) and satisfies a minimum-like (resp. maximum-like) principle. This is summarised in the following theorem.
\begin{thm}\label{thm1}
	Under the assumptions in \cite{DeDona_Levine2013barriers} and \cite{Ester_Asch_Streif_2019}, every integral curve $x^{\bar{u}}$ on $\DAM$ (resp. $\DMM$) and the corresponding input function $\bar{u}$ satisfy the following necessary conditions. There exists a nonzero absolutely continuous maximal solution $\lambda^{\bar{u}}$ to the adjoint equation:
	\begin{align}
	&\dot{\lambda}^{\bar{u}}(t) = -\left( \frac{\partial f}{\partial x}(x^{\bar{u}}(t),\bar{u}(t)) \right)^T \lambda^{\bar{u}}(t),\nonumber
	\end{align}
	such that
	\begin{align}		
	&\min_{u\in U}\{\lambda^{\bar{u}}(t)^Tf(x^{\bar{u}}(t),u) \} &&= \lambda^{\bar{u}}(t)^T f(x^{\bar{u}}(t),\bar{u}(t)) \nonumber \\
	&															 &&= 0 \label{eq_Hamil_A} \\
	%\end{align}
	%\begin{align}	
	\Big(\Big.\text{resp.}\;
	&\max_{u\in U}\{\lambda^{\bar{u}}(t)^Tf(x^{\bar{u}}(t),u) \} &&= \lambda^{\bar{u}}(t)^T f(x^{\bar{u}}(t),\bar{u}(t)) \nonumber \\ 
	&															&&= 0\label{eq_Ham_max_M}
	\Big.\Big),
	\end{align}
	for almost all $t$. Moreover, if $x^{\bar{u}}$ intersects $G_0$ in finite time, we have:
	%\begin{equation}
	$\lambda^{\bar{u}}(\bar{t}) = (Dg_{i^*}(z))^T,$ %\label{eq_adj}
	%\end{equation}
	where
	\begin{align}
	&\min_{u\in U}\max_{i\in\mathbb{I}(z)} L_fg_i(z,u) = L_fg_{i^*}(z,\bar{u}(\bar{t})) = 0 \label{thm1_ult_tan_A} \\
	\Big(\Big.\text{resp.}\quad
	&\max_{u\in U}\max_{i\in\mathbb{I}(z)} L_fg_i(z,u) = L_fg_{i^*}(z,\bar{u}(\bar{t})) = 0 \label{thm1_ult_tan_M}
	\Big.\Big),
	\end{align}
	$\bar{t}$ denotes the time at which $G_0$ is reached, and $z \triangleq x^{(\bar{u},\bar{x},t_0)}(\bar{t})\in G_0$.
\end{thm}
In addition to using the necessary conditions to analyse the predator-prey model (for which we will exactly describe the special control $\bar{u}$, and present conditions under which certain parts of the set $\A$ exists), we will also use them to construct the sets as follows:
first, identify points of \emph{ultimate tangentiality} on $G_0$, via \eqref{thm1_ult_tan_A} (resp. \eqref{thm1_ult_tan_M}).
Second, determine the input realisation associated with $\DAM$ (resp. $\DMM$) using the \emph{Hamiltonian minimisation} (resp. \emph{maximisation}) condition \eqref{eq_Hamil_A} (resp. \eqref{eq_Ham_max_M}).
Third, integrate the system dynamics and the adjoint equations backwards in time from the points of ultimate tangentiality to obtain the barrier curves.
%\begin{itemize}
%	\item Identify points of \emph{ultimate tangentiality} located on $G_0$, via the conditions \eqref{thm1_ult_tan_A} (resp. \eqref{thm1_ult_tan_M}). 
%	\item Determine the input realisation associated with $\DAM$ (resp. $\DMM$) using the \emph{Hamiltonian minimisation} (resp. \emph{maximisation}) condition \eqref{eq_Hamil_A} (resp. \eqref{eq_Ham_max_M}).
%	\item Obtain the barrier curves by integrating the system dynamics and the adjoint equations backwards in time from the points of ultimate tangentiality.
%\end{itemize}

Since the conditions in Theorem \ref{thm1} are necessary, it may happen that some integral curves obtained through this method (or parts of them) do not define parts of the boundary of $\A$ or $\M$ and may have to be ignored. Thus, in general we refer to trajectories obtained from these conditions as \emph{candidate barrier} and \emph{candidate invariance barrier} trajectories.

\section{Set-Based Analysis of the Lotka-Volterra Equations} \label{sec:3}
In this section we present the analysis of the aforementioned sets. We also derive conditions under which candidate barrier trajectories exist.

Consider the Lotka-Volterra equations \eqref{eq_LotkaVolterra},
the system's axes form two trivial robustly invariant manifolds, with the direction of flow along them determined by the signs of $\alpha + u_1^{\max}$ and $-\gamma + u_2^{\max}$. We use the necessary conditions of Theorem~\ref{thm1} to show that the constrained Lotka-Volterra model has a non-trivial MRPI if and only if $u_i^{\min} = u_i^{\max}$, $i=1,2$. We use the well-known fact that with a constant input the integral curves of the Lotka-Volterra model trace out periodic orbits, for example, see Lemma 1 in \cite{bayen2019minimalTimeCrisisToViabKernel}.
\begin{prop}
	The Lotka-Volterra model as in \eqref{eq_LotkaVolterra} with $x(t)\in G\subset[0,\infty[\times[0,\infty[$ for all $t$, with $G$ as defined in \eqref{def_G}, and $u(t)\in U\triangleq[u_1^{\min}, u_1^{\max}]\times [u_2^{\min}, u_2^{\max}]$, has an MRPI that includes points not on the axes, if and only if $U$ is a singleton.
\end{prop}
\begin{pf}
	Suppose the system has an MRPI as described, labelled $\M$, and let us concentrate on the invariance barrier trajectories on $\DMM$. From Theorem~\ref{thm1}, any curve on $\DMM$ satisfies the Hamiltonian maximisation condition \eqref{eq_Ham_max_M} for almost all $t$. Suppose now that there exists a $\tilde{u}\in U$ such that $\lambda^{\bar{u}}(\hat{t})^Tf(x^{\bar{u}}(\hat{t}),\tilde{u}) < 0$ where $\hat{t}$ is an arbitrary Lebesgue point of $\bar{u}$. Then, we could specify the constant input $u(t)\equiv \arg\min \lambda^{\bar{u}}(\hat{t})^Tf(x^{\bar{u}}(\hat{t}),u)$, which would result in a periodic orbit that intersects $\Mcomp$, contradicting the fact that $\M$ is robustly invariant. We conclude that for any trajectory on $\DMM$, we have $\lambda^{\bar{u}}(t)^Tf(x^{\bar{u}}(t),u) = 0$ for all $u\in U$ for almost all $t$. In particular,
	\begin{equation}
		\max_{u\in U}\lambda^{\bar{u}}(t)^Tf(x^{\bar{u}}(t),u) = \min_{u\in U}\lambda^{\bar{u}}(t)^Tf(x^{\bar{u}}(t),u) = 0,\label{eq_max_equals_min}
	\end{equation}
	for almost all $t$.
	
	\sloppy
	Let $\hat{u}(t)=(\hat{u}_1(t), \hat{u}_2(t))^T\triangleq \arg\max_{u} \lambda^{\bar{u}}(t)^Tf(x^{\bar{u}}(t),u)$ and $\check{u}(t)=(\check{u}_1(t), \check{u}_2(t))^T\triangleq\arg\min_{u} \lambda^{\bar{u}}(t)^Tf(x^{\bar{u}}(t),u)$. From the form of the Hamiltonian we see that the functions $\hat{u}_i$ and $\check{u}_i$ are always saturated, that is $\hat{u}_i(t)\in\{u_i^{\min},u_i^{\max}\}$ and $\check{u}_i(t)\in\{u_i^{\min},u_i^{\max}\}$, $i=1,2$. From \eqref{eq_max_equals_min} we can conclude that:
	\[
		\lambda_1^{\bar{u}}(t)x_1^{\bar{u}}(t)(u^{\min}_1- u^{\max}_1) + \lambda_2^{\bar{u}}(t) x_2^{\bar{u}} (t)(u^{\min}_2 - u^{\max}_2) = 0,
	\]
	for almost all $t$. We now argue that $\lambda_1^{\bar{u}}(t)$ and $\lambda_2^{\bar{u}}(t)$ are nonzero for almost all time. Indeed, suppose $\lambda_1^{\bar{u}}(t) = 0$ over an interval of time. Then, $\dot{\lambda}^{\bar{u}}_1(t) = 0$ over this same interval, which gives $-\delta x_2(t) \lambda_2^{\bar{u}}(t) = 0$, and thus, because the state is assumed to be positive, $\lambda_2^{\bar{u}}(t) = 0$ over this same interval. But this is impossible because $\lambda^{\bar{u}}(t) \neq 0$, from Theorem~\ref{thm1}. A similar argument holds for $\lambda_2^{\bar{u}}$. Therefore, \eqref{eq_max_equals_min} holds if and only if $u_1^{\max} = u_1^{\min}$ and $u_2^{\max} = u_2^{\min}$ for almost all $t$.
	
	Turning our attention to the set $\DMO$, by a similar argument we require that
	\begin{equation}
		\min_{u\in U} \max_i L_f g_i(z,u) = \max_{u\in U} \max_i L_f g_i(z,u) = 0\label{eq_points_on_DMO}
	\end{equation}
	for points on $\DMO$. Thus, either i) $Dg_{i^*1} = 0$ and $u_2^{\max} = u_2^{\min}$; or ii) $Dg_{i^*2} = 0$ and $u_1^{\max} = u_1^{\min}$; or iii) $u_1^{\max} = u_1^{\min}$ and $u_2^{\max} = u_2^{\min}$, where $Dg_{i^*} = (Dg_{i^*1}, Dg_{i^*2})$, must hold on $\DMO$. Focusing on case i), we see that this condition only holds if $x_1 = \frac{\gamma - u_2}{\delta}$, (with $u_2 = u_2^{\max} = u_2^{\min}$) and thus the condition only holds at an isolated point on the active constraint $g_{i^*}$. In order for this point to be robustly invariant, we would require $\dot{x}_1 = 0$, which can only be true if $u_1^{\max} = u_1^{\min}$. A similar argument for case ii) leads to $u_2^{\max} = u_2^{\min}$ for points on $\DMO$. This completes the proof.
\end{pf}

\fussy
Though the constrained Lotka-Volterra model as specified above cannot have a nontrivial MRPI, it can have an admissible set, as we will show. We continue our analysis with the same constraints on $u$, $u_i(t)\in[u_i^{\min},u_i^{\max}]$, and we let the state be constrained to a box: $x_i(t) \in [\underline{x}_i, \overline{x}_i]$ for all $t$, with $ 0 < \underline{x}_i < \overline{x}_i$.

\subsection{Points of ultimate tangentiality for the admissible set}
We label the state constraints as follows: $g_1(x) = x_1 - \overline{x}_1$, $g_2(x) = - x_1 + \underline{x}_1$, $g_3(x) = x_2 - \overline{x}_2$, $g_4(x) = - x_2 + \underline{x}_2$, and the $i$-th ultimate tangentiality point $z_i\triangleq (z_{i1}, z_{i2})^T\in\mathbb{R}^2$. First, we consider the constraints for the prey, i.e. $g_1(x)$ and $g_2(x)$. Invoking condition \eqref{thm1_ult_tan_A} we get:
\begin{align*}
&\min_{u\in U } L_f g_1(z,u) =
 \min_{u\in U } \left\{ \alpha\overline{x}_1-\beta\overline{x}_1 z_{12} + u_1 \overline{x}_1 \right\} = 0, \\
&\min_{u\in U } L_f g_2(z,u)  =
 \min_{u\in U } \left\{ -\alpha\underline{x}_1 + \beta\underline{x}_1 z_{22} - u_1 \underline{x}_1  \right\} = 0,
\end{align*}
from where we get the two points $z_1 = (\overline{x}_1, \frac{\alpha + u_1^{\min}}{\beta})$ for $g_1$ and  $z_2 = (\underline{x}_1, \frac{\alpha + u_1^{\max}}{\beta})$ for $g_2$. Now we focus on the predator constraints, $g_3(x)$ and $g_4(x)$. Again invoking condition \eqref{thm1_ult_tan_A} we get:
\begin{align*}
&\min_{u\in U } L_f g_3(z,u)  =
\min_{u\in U } \left\{ \delta z_{31} \overline{x}_2 - \gamma \overline{x}_2 + u_2 \overline{x}_2 \right\} = 0, \\
&\min_{u\in U }L_f g_4(z,u)  =
\min_{u\in U } \left\{ -\delta z_{41} \underline{x}_2 + \gamma \underline{x}_2 - u_2 \underline{x}_2 \right\} = 0,
\end{align*}
from where we identify $z_3 = ( \frac{\gamma-u_2^{\min}}{\delta}, \overline{x}_2)$ for $g_3$ and  $z_4 = (\frac{\gamma-u_2^{\max}}{\delta}, \underline{x}_2)$ for $g_4$.

\subsection{Input realisation associated with the barrier}
Invoking condition \eqref{eq_Hamil_A} to determine the input associated with the barrier, we get:
\begin{align*}
	\min_{u\in U}\left\{\lambda(t)^Tf(x(t),u) \right\} 
	&= \min_{u\in U} \{ \lambda_1 \left( \alpha x_1 - \beta x_1 x_2 + u_1x_1 \right) \\ 
	&+ \lambda_2 \left( \delta x_1 x_2 - \gamma x_2 +  u_2x_2 \right) \} = 0 
\end{align*}
for almost every $t$, with $\lambda \triangleq (\lambda_1, \lambda_2)^T$. We get:
\begin{align}
\bar{u}_1(t) = 
\begin{split}
\begin{cases}
\label{input_real_u1_A}
u_1^{\min} & \text{if} \quad \lambda_1(t) \geq 0\\
u_1^{\max} & \text{if} \quad \lambda_1(t) < 0,
\end{cases}
\end{split}
\end{align}
and
\begin{align}
\bar{u}_2(t) = 
\begin{split}
\begin{cases}
\label{input_real_u2_A}
u_2^{\min} & \text{if} \quad \lambda_2(t) \geq 0\\
u_2^{\max} & \text{if} \quad \lambda_2(t) < 0.
\end{cases}
\end{split}
\end{align}
The adjoint equation is
\begin{align} \label{eq_adj_lotka_volterra}
\dot{\lambda}
=
\begin{pmatrix}
-\alpha+\beta x_2 - u_1 & -\delta x_2\\
\beta x_1 & -\delta x_1 + \gamma - u_2
\end{pmatrix}
\lambda,
\end{align}
with the final conditions (from $\lambda(\bar{t}) = (Dg_{i^*}(z))^T),$  $\lambda(\bar{t}) = (1,0)^T$, $(-1,0)^T$, $(0,1)^T$, $(0,-1)^T$ associated with $z_1$, $z_2$, $z_3$ and $z_4$, respectively.

We note that there are four lines in the state space (that intersect the four points of ultimate tangentiality) where the control $\bar{u}$ switches, summarised in the following Proposition. 

\begin{prop}\label{prop_3_switches}
	Switches in $\bar{u}$ occur as follows:
	\begin{itemize}
		\item If $x^{\bar{u}}(t)\in \{(x_1,x_2) : \delta x_1 - \gamma + u_2^{\max} = 0\}$  and $\lambda_2(t) < 0$, then $\bar{u}_1$ switches from $u_1^{\max}$ to $u_1^{\min}$.
		\item If $x^{\bar{u}}(t)\in \{(x_1,x_2) : \delta x_1 - \gamma + u_2^{\min} = 0\}$  and $\lambda_2(t) > 0$, then $\bar{u}_1$ switches from $u_1^{\min}$ to $u_1^{\max}$.
		\item If $x^{\bar{u}}(t)\in \{(x_1,x_2) : \alpha -\beta x_2  + u_1^{\max} = 0\}$  and $\lambda_1(t) < 0$, then $\bar{u}_2$ switches from $u_2^{\min}$ to $u_2^{\max}$.
		\item If $x^{\bar{u}}(t)\in \{(x_1,x_2) : \alpha -\beta x_2 + u_1^{\min} = 0\}$  and $\lambda_1(t) > 0$, then $\bar{u}_2$ switches from $u_2^{\max}$ to $u_2^{\min}$.
	\end{itemize}
	
\end{prop}
\begin{pf}
	From condition \eqref{eq_Hamil_A}, we see that if $\lambda_1(t) = 0$, then $\lambda_2(t) (\delta x_1(t) x_2(t) - \gamma x_2(t) + \bar{u}_2(t) x_2(t)) = 0$, and using the fact that $\lambda_2(t) \neq 0$ and $x_2(t) \neq 0$, we see that if $\lambda_2(t) < 0$ then a switch in $\bar{u}_1(t)$ occurs on the line segment given by $\{(x_1,x_2) : \delta x_1 - \gamma + u_2^{\max}\}$. Because $\lambda_2(t) < 0$, we have $\dot{\lambda}_1(t) > 0$, implying $\lambda_1(t) < 0$ on an interval before $t$, implying $\bar{u}_1(t) = u_1^{\max}$ before $t$. This same argument carries over to the remaining three cases.
\end{pf}

\subsection{Existence of candidate barrier trajectories}
As already mentioned in Section \ref{sec:2} the conditions stated in Theorem \ref{thm1} are necessary. Therefore, some integral curves or parts of them may need to be ignored since they do not define the boundary of the sets. In particular, integral curves evolving outside the constrained state space need to be ignored. The next proposition gives conditions under which a point of ultimate tangentiality is indeed associated with a candidate barrier trajectory evolving backwards into $G_-$.
\begin{prop}\label{prop_existence}
	There exists a candidate barrier trajectory associated with $\A$, partly contained in $G_-$ and ending at the point of ultimate tangentiality
%	\begin{itemize}
%		\item $z_1 \Leftrightarrow \left( \delta \overline{x}_1 - \gamma + u_2^{\max}\right) > 0$,
%		\item $z_2 \Leftrightarrow \left( \delta \underline{x}_1 - \gamma + u_2^{\min}\right) < 0$,
%		\item $z_3 \Leftrightarrow \left( \alpha - \beta\overline{x}_2 + u_1^{\min}\right) < 0$,
%		\item $z_4 \Leftrightarrow \left( \alpha - \beta\underline{x}_2 + u_1^{\max}\right) > 0$.
%	\end{itemize}
\begin{align*}
&z_1 \Leftrightarrow \left( \delta \overline{x}_1 - \gamma + u_2^{\max}\right) > 0, \; 
z_2 \Leftrightarrow \left( \delta \underline{x}_1 - \gamma + u_2^{\min}\right) < 0, \\
&z_3 \Leftrightarrow \left( \alpha - \beta\overline{x}_2 + u_1^{\min}\right) < 0, \; 
z_4 \Leftrightarrow \left( \alpha - \beta\underline{x}_2 + u_1^{\max}\right) > 0.
\end{align*}
\end{prop} 
\begin{pf}
	Consider $g_1$ and the corresponding point of ultimate tangentiality $z_1$ along with the final adjoint $\lambda(\bar{t}) = (1,0)^T$.
	From \eqref{eq_adj_lotka_volterra} it follows that $\dot{\lambda}_2(\bar{t}) = \beta \overline{x}_1\lambda_1(\bar{t}) > 0$, which implies that $\lambda_2 (t) < 0$ and $\bar{u}_2(t) = u_2^{\max}$ over a time interval before $\bar{t}$. Because $\lambda_2(t)<0$ over this interval, the integral curve associated with $\bar{u}$ will evolve backwards into $G_-$ if and only if $\dot{x}_2(\bar{t}) > 0$, which gives us the first statement. Similar arguments carry over to the other three statements, which completes the proof.
\end{pf} 
\begin{figure}[h]
	\begin{center}
		\includegraphics[width=4\linewidth,height=5cm,keepaspectratio]{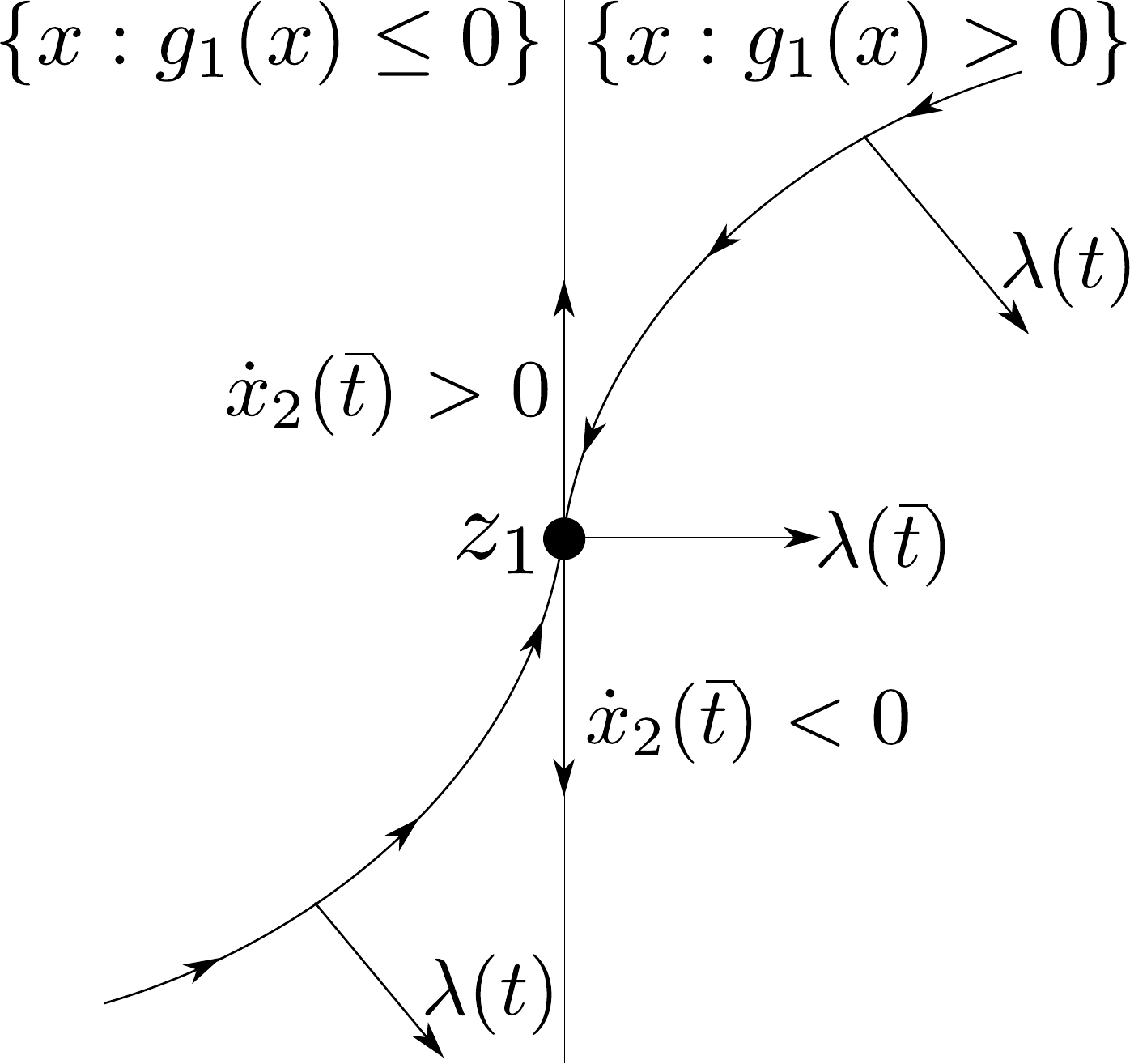} 
		\caption{Clarification of the proof of Proposition~\ref{prop_existence}: knowing that $\lambda_2(t)<0$ before the barrier intersects $z_1$, we must have $\dot{x}_2(\bar{t}) >0$.} 
		\label{fig:prop_3}                             
	\end{center}                        
\end{figure}

\section{Numerical Example} \label{sec_examples}
\begin{figure}[h]
	\begin{center}\
		\includegraphics[width=7.5\linewidth,height=7.5cm,keepaspectratio]{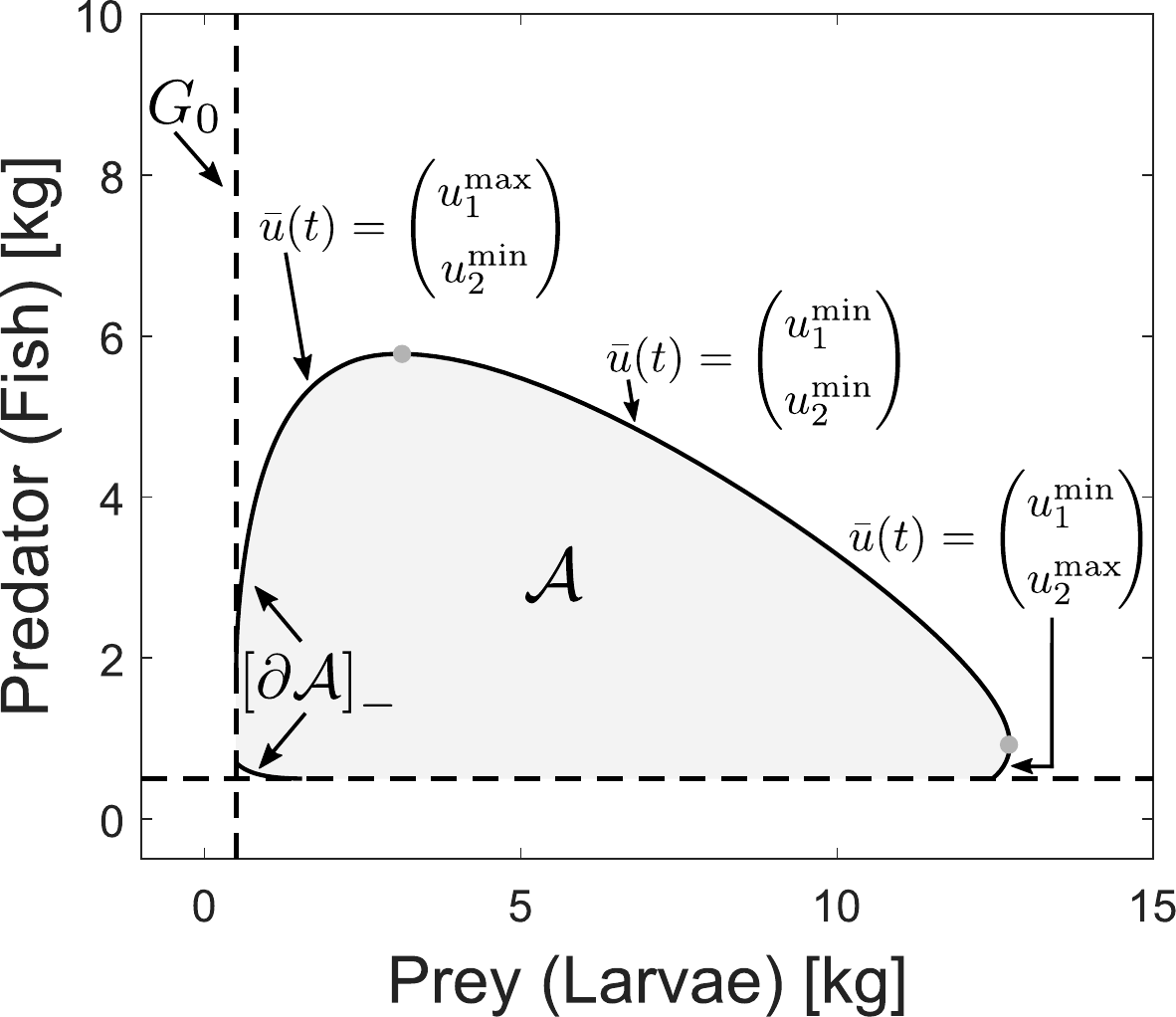} 
		\caption{ 
			Admissible set for a system with interacting larvae and fish to prevent underpopulation. 
			%The input $\bar{u}(t)$ associated with the trajectory initiating from $z_2$ switches at the grey points according to Proposition \ref{prop_3_switches}. 
			The grey points indicate the switching of $\bar{u}(t)$, associated with the barrier trajectory initiating from $z_2$.}%, according to Proposition \ref{prop_3_switches}} 
		\label{fig:Example_larvae_fish_underpop}                             
	\end{center}                              
\end{figure}

This example is based on the recent trend in aquaculture to feed fish with insect larvae-based food for efficient fish production. The environments of the species are often decoupled. Nevertheless, we assume that the fish are continuously supplied with larvae on which to feed. 
%From the predator-prey model \eqref{eq_LotkaVolterra}, in this example, the prey are represented by the larvae of \textit{Hermetia Illucens} insects and the predators by \textit{Oreochromis niloticus} fish. 
The prey are represented by the larvae of \textit{Hermetia illucens} insects and the predators by \textit{Oreochromis niloticus} fish. 
The goal of this example is to show the application of the proposed set-based method for sustainable production of larvae and fish.
First, we impose state constraints to prevent underpopulation of both species with $g_2(x) \triangleq -x_1 + 0.5$ and $g_4(x) \triangleq -x_2 + 0.5$. 
The input bounds are $u_1 \in [10,20]$ and $u_2 \in [-10,-20]$. 
%We impose the state constraints $g_1(x) \triangleq x_1 - 10$, $g_2(x) \triangleq -0.5 + x_1$, $g_3(x) \triangleq x_2 - 10$, $g_4(x) \triangleq -0.5 + x_2$ as well as the input bounds $u_1 \in [10,20]$ and $u_2 \in [-10,-20]$.
%The parameters for the model are given in the appendix with a detailed explanation.
The model parameters $\alpha, \beta, \gamma$ and $\delta$ are derived from the parameters $k_1$, $k_4$ \citep{diener2009}, $k_2$, $k_3$, $k_5$ \citep{terpstra2015} and $\eta_1$ \citep{rana2015,stamer2014black} as listed in Table \ref{tb:fish_larvae_raw}.
\begin{table}[!h]
	\caption{Lotka-Volterra model parameters for fish (\textit{Oreochromis niloticus}) and larvae (\textit{Hermetia illucens}) production system}
	\begin{tabular}{|l|l|l|l|}
		\hline
		\multicolumn{2}{|l|}{parameter and meaning} & value   & unit  \\ \hline
		$k_1$    & larval growth from egg to adult      & 30            & $[d]$          \\ \hline
		$k_2$    & fish growth from egg to adult        & 225           & $[d]$          \\ \hline
		$k_3$    & larvae consumption rate of fish      & 22            & $[gd^{-1}]$    \\ \hline
		$k_4$    & weight of  adult larva               & 0.120         & $[g]$          \\ \hline
		$k_5$    & weight of adult fish                 & 700           & $[g]$          \\ \hline
		$\eta_1$ & conversion eff. of larval food       & 0.55          & [-]            \\ \hline
		$\alpha$ & $1/k_1$                              & 0.0014        & $[h^{-1}]$     \\ \hline
		$\beta$  & $k_3/(k_4k_5)$                       & 10.8          & $[{kgh}^{-1}]$ \\ \hline
		$\delta$ & $\eta_1\beta$                        & 5.94          & $[{kgh}^{-1}]$ \\ \hline
		$\gamma$ & $1/k_2$                              & \num{1.85e-5} & $[h^{-1}]$     \\ \hline
	\end{tabular}
	\label{tb:fish_larvae_raw}	
\end{table}
% USAGE OF PARAMETER PER HOUR beta = 0.003 * 3600 = 10.8 and delta = 0.0018 * 3600 = 6.48
Since $\alpha$ and $\gamma$ are vastly smaller than $\beta$ and $\delta$ as well as the input bounds, they do not play a decisive role for the dynamics. Hence, we set them to zero.
%We identify: $z_1 \approx (10, 0.926)^T$, $z_2 \approx (0.5, 1.852)^T$, $z_3 \approx (3.086, 10)^T$, and $z_4 \approx (1.543, 0.5)^T$.
We identify:  $z_2 \approx (0.5, 1.852)^T$ and $z_4 \approx (1.543, 0.5)^T$.
By integrating backwards, we find two candidate barrier trajectories defining the boundary of $\A$, shown in Figure \ref{fig:Example_larvae_fish_underpop}.

\begin{figure}[h]
	\begin{center}\
		\includegraphics[width=6.5\linewidth,height=7.5cm,keepaspectratio]{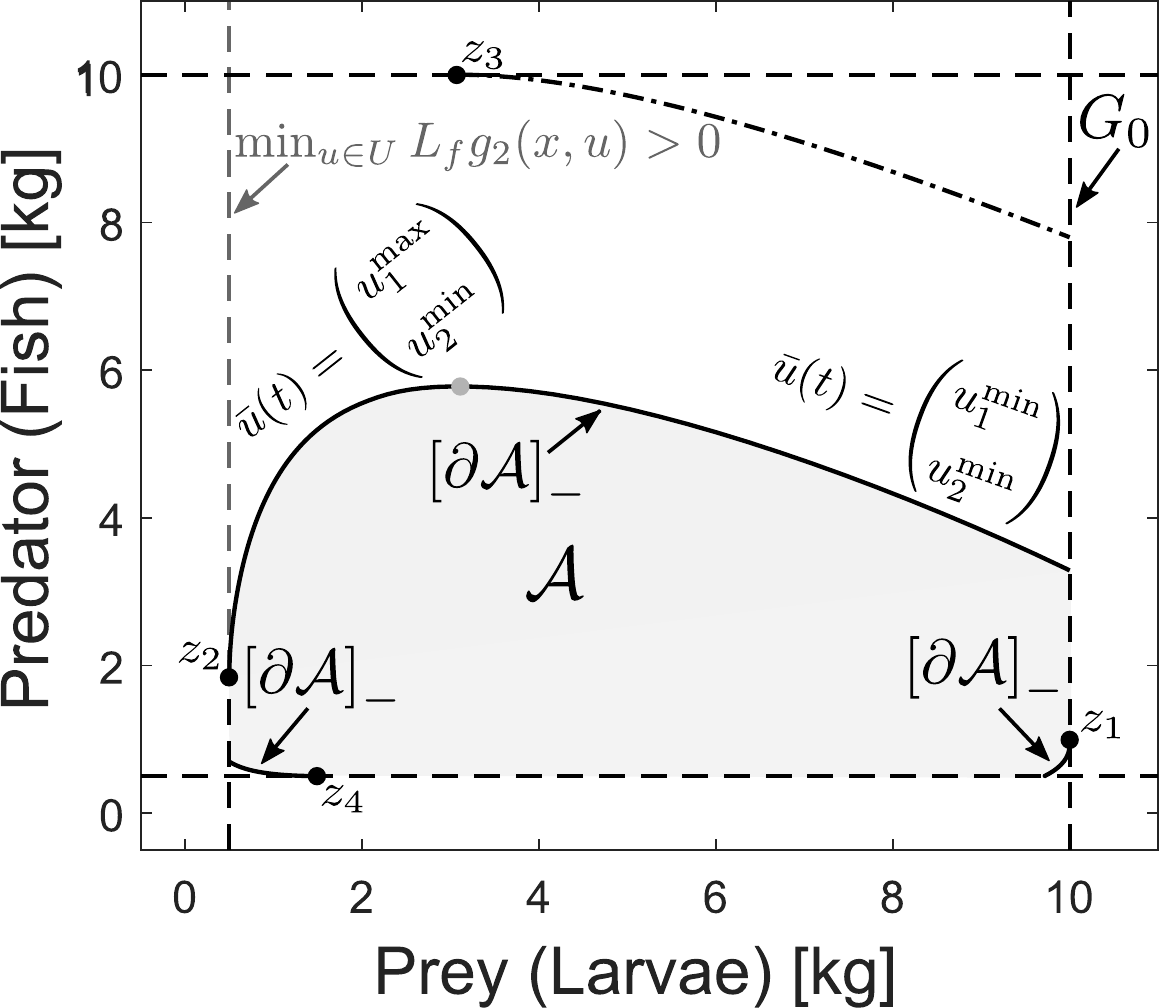} 
		\caption{Admissible set for a system with interacting larvae and fish. The input realisation $\bar{u}(t)$ associated with the barrier trajectory initiating from $z_2$ switches at the grey point.} 
		\label{fig:Example_larvae_fish}                             
	\end{center}                              
\end{figure}

Next, we additionally impose state constraints to prevent overpopulation with $g_1(x) \triangleq x_1 - 10$ and $g_3(x) \triangleq x_2 - 10$.
We identify: $z_1 \approx (10, 0.926)^T$ and $z_3 \approx (3.086, 10)^T$ and integrate backwards from all four points of tangentiality.
We find four candidate barrier trajectories, shown in Figure \ref{fig:Example_larvae_fish}. 
The difference to the example before is the existence of a trajectory evolving into $G_-$ which does not define a part of the boundary of the admissible set.
Hence, we have to ignore the dash-dotted line drawn from $z_3$ and ending on the set $\{x:g_1(x) = 0\}$ because if we were to take this curve as the set's boundary, we would include points on the set $\{x:g_2(x) = 0\}$ for which $\min_{u\in U} L_fg_2(x,u) > 0$ as being on the admissible set's boundary, which would be impossible. Thus, the boundary of the admissible set is defined by the three integral curves ending on $z_1$, $z_2$, and $z_4$ shown by the continuous lines.

\section{Discussion}\label{Sec_discussion}
% Simplicity
It is interesting to note that other works in the literature that concentrate on biological applications of the viability kernel describe its boundary, and argue that it is made up of special integral curves of the system. A similar statement is made in the theory of barriers, the main observation being that these curves satisfy a minimum/maximum-like principle, as described in Theorem~\ref{thm1}. 
In this new paradigm it can be seen that the description of the boundary is much simpler than the arguments in, for example, \cite{bayen2019minimalTimeCrisisToViabKernel}.

%Discussion of sets for predator-prey system
Recall that the considered predator-prey system never has a (nontrivial) MRPI. Thus, the interpretation is that if one can control the introduction or harvesting of species, then this has to be done carefully, because there always exists an input (intervention strategy) that results in over- or under-population of a species in the future. 
Thus, it is vital that the population number remains in $\A$. 
%Otherwise, if the state enters $\Acomp$, there is \emph{no means} of preventing over- or under-population. Lastly, if the state is located on $\DAM$, then there exists a \emph{unique} input that prevents a constraint violation in the future. 
Furthermore, the unique input on the boundary of $\A$ is always saturated, and it switches according to the statement in Proposition~\ref{prop_3_switches}. %\textcolor{red}{Probably the meaning of this statement in terms of the example ??}
%Regarding to the example, on the boundary of $\A$ only one intervention strategy will lead to not result in an over- or underpopulation in the future and this strategy changes at the switching points.
For a given system, similar to the model of the CEA system presented in this work, the proposed method can be used to set up the initial population as well as intervention strategies through the addition and/or removal of the species.

\section{Conclusion}\label{sec_conclusion}
%We used the theory of barriers to analyse the admissible set (also known as viability kernel) and the MRPI of predator-prey systems modelled by the Lotka-Volterra equations. 
%We have reproduced a result from the literature, emphasising the simplicity of our approach, and discovered new aspects, such as the nonexistence of nontrivial MRPIs in the constrained Lotka-Volterra model. 
%We also present an example for an ecosystem under controlled conditions for food production. 
%In the future it could be interesting to analyse more biological systems using these techniques.

The main goal of this paper was to analyse CEA systems with interacting species, modelled by the Lotka-Volterra equations, to identify the admissible set and the MRPI which yields information about sustainable intervention strategies. 
We used the theory of barriers to determine special trajectories that define the boundaries of the admissible set and discovered new aspects, such as the nonexistence of nontrivial MRPIs in the considered constrained Lotka-Volterra model. 
Furthermore, we obtained conditions of the systems' parameters that guarantees the existence of these special trajectories.
We illustrated our results based on a CEA system with larve and fish as interacting species.
Future research could focus on extending the results to robust admissible sets, as in \cite{regnier2015robust}, in the context of viability theory. Another interesting idea could be to analyse the sets of an extended Lotka-Volterra model that describes the effects of humans feeding one species to another.

%\begin{ack}
%	\vspace{-2mm}
%This project has received funding from the European Social Fund (ESF).
              	%\vspace{-2mm}  
%\end{ack}

\bibliography{bib/bib_predprey}             % bib file to produce the bibliography

\begin{thebibliography}{20}
\providecommand{\natexlab}[1]{#1}
\providecommand{\url}[1]{\texttt{#1}}
\providecommand{\urlprefix}{URL }
\expandafter\ifx\csname urlstyle\endcsname\relax
  \providecommand{\doi}[1]{doi:\discretionary{}{}{}#1}\else
  \providecommand{\doi}{doi:\discretionary{}{}{}\begingroup
  \urlstyle{rm}\Url}\fi

\bibitem[{Barnosky et~al.(2011)Barnosky, Matzke, Tomiya, Wogan, Swartz,
  Quental, Marshall, McGuire, Lindsey, Maguire, Mersey, and
  Ferrer}]{Barnosky2011}
Barnosky, A.D., Matzke, N., Tomiya, S., Wogan, G.O.U., Swartz, B., Quental,
  T.B., Marshall, C., McGuire, J.L., Lindsey, E.L., Maguire, K.C., Mersey, B.,
  and Ferrer, E.A. (2011).
\newblock Has the earth’s sixth mass extinction already arrived?
\newblock \emph{Nature}, 471(7336), 51--57.

\bibitem[{Bayen and Rapaport(2019)}]{bayen2019minimalTimeCrisisToViabKernel}
Bayen, T. and Rapaport, A. (2019).
\newblock Minimal time crisis versus minimum time to reach a viability kernel:
  A case study in the prey-predator model.
\newblock \emph{Optimal Control Applications and Methods}, 40(2), 330--350.

\bibitem[{Conrad et~al.(2017)Conrad, Daniel, and Vincent}]{CONRAD2017}
Conrad, Z., Daniel, S., and Vincent, V. (2017).
\newblock Vertical farm 2.0: Designing an economically feasible vertical farm -
  a combined european endeavor for sustainable urban agriculture.
\newblock Technical report, Association for Vertical Farming.
\newblock White Paper.

\bibitem[{De~Dona and L{\'e}vine(2013)}]{DeDona_Levine2013barriers}
De~Dona, J.A. and L{\'e}vine, J. (2013).
\newblock On barriers in state and input constrained nonlinear systems.
\newblock \emph{SIAM Journal on Control and Optimization}, 51(4), 3208--3234.

\bibitem[{De~Lara and Doyen(2008)}]{deLara2008sustainableBOOK}
De~Lara, M. and Doyen, L. (2008).
\newblock \emph{Sustainable management of natural resources: mathematical
  models and methods}.
\newblock Springer Science \& Business Media.

\bibitem[{Diener et~al.(2009)Diener, Zurbr{\"u}gg, and Tockner}]{diener2009}
Diener, S., Zurbr{\"u}gg, C., and Tockner, K. (2009).
\newblock Conversion of organic material by black soldier fly larvae:
  establishing optimal feeding rates.
\newblock \emph{Waste Management \& Research}, 27(6), 603--610.

\bibitem[{Eisenack et~al.(2006)Eisenack, Scheffran, and
  Kropp}]{eisenack2006viabilityAnalysis}
Eisenack, K., Scheffran, J., and Kropp, J.P. (2006).
\newblock Viability analysis of management frameworks for fisheries.
\newblock \emph{Environmental Modeling \& Assessment}, 11(1), 69--79.

\bibitem[{Esterhuizen et~al.(2019)Esterhuizen, Aschenbruck, and
  Streif}]{Ester_Asch_Streif_2019}
Esterhuizen, W., Aschenbruck, T., and Streif, S. (2019).
\newblock On maximal robust positively invariant sets in constrained nonlinear
  systems.
\newblock \emph{arXiv:1904.01985}.

\bibitem[{Fort et~al.(2017)Fort, Dieguez, Halty, and Lima}]{FORT201723}
Fort, H., Dieguez, F., Halty, V., and Lima, J.M.S. (2017).
\newblock Two examples of application of ecological modeling to agricultural
  production: Extensive livestock farming and overyielding in grassland
  mixtures.
\newblock \emph{Ecological Modelling}, 357, 23--34.

\bibitem[{Funabashi(2018)}]{Funabashi2018}
Funabashi, M. (2018).
\newblock Human augmentation of ecosystems: objectives for food production and
  science by 2045.
\newblock \emph{npj Science of Food}, 2(1), 1--11.

\bibitem[{Getz(2012)}]{getz2012biomass}
Getz, W.M. (2012).
\newblock A biomass flow approach to population models and food webs.
\newblock \emph{Natural resource modeling}, 25(1), 93--121.

\bibitem[{Gligorescu et~al.(2019)Gligorescu, Toft, Hauggaard-Nielsen, Axelsen,
  and Nielsen}]{Anton2019}
Gligorescu, A., Toft, S., Hauggaard-Nielsen, H., Axelsen, J.A., and Nielsen,
  S.A. (2019).
\newblock Development, growth and metabolic rate of hermetia illucens larvae.
\newblock \emph{Journal of Applied Entomology}, 143(8), 875--881.

\bibitem[{Graham et~al.(2007)Graham, Knapp, Van~Vleck, Bloor, Lane, and
  Graham}]{Graham2007}
Graham, D.W., Knapp, C.W., Van~Vleck, E.S., Bloor, K., Lane, T.B., and Graham,
  C. (2007).
\newblock Experimental demonstration of chaotic instability in biological
  nitrification.
\newblock 385--393.

\bibitem[{Martinet and Doyen(2007)}]{martinet2007sustainabilityViableControl}
Martinet, V. and Doyen, L. (2007).
\newblock Sustainability of an economy with an exhaustible resource: A viable
  control approach.
\newblock \emph{Resource and energy economics}, 29(1), 17--39.

\bibitem[{Rana et~al.(2015)Rana, Salam, Hashem, and Islam}]{rana2015}
Rana, K.S., Salam, M., Hashem, S., and Islam, M.A. (2015).
\newblock Development of black soldier fly larvae production technique as an
  alternate fish feed.
\newblock \emph{International Journal of Research in Fisheries and
  Aquaculture}, 5(1), 41--47.

\bibitem[{Regnier and De~Lara(2015)}]{regnier2015robust}
Regnier, E. and De~Lara, M. (2015).
\newblock Robust viable analysis of a harvested ecosystem model.
\newblock \emph{Environmental Modeling \& Assessment}, 20(6), 687--698.

\bibitem[{Schmidt et~al.(2014)Schmidt, Barney, Williams, Bessin, Coolong, and
  Harwood}]{Schmidt2014}
Schmidt, J.M., Barney, S.K., Williams, M.A., Bessin, R.T., Coolong, T.W., and
  Harwood, J.D. (2014).
\newblock Predator–prey trophic relationships in response to organic
  management practices.
\newblock \emph{Molecular Ecology}, 23(15), 3777--3789.

\bibitem[{Stamer et~al.(2014)Stamer, Wessels, Neidigk, and
  Hoerstgen-Schwark}]{stamer2014black}
Stamer, A., Wessels, S., Neidigk, R., and Hoerstgen-Schwark, G. (2014).
\newblock Black soldier fly (hermetia illucens) larvae-meal as an example for a
  new feed ingredients’ class in aquaculture diets.
\newblock \emph{Proceedings of the 4th ISOFAR Scientific Conference}.

\bibitem[{Terpstra(2015)}]{terpstra2015}
Terpstra, A.H. (2015).
\newblock Feeding and growth parameters of the tilapia (oreochromis niloticus)
  in the body weight range from newly hatched larvae (about 5 mgrams) to about
  700 grams.
\newblock \emph{Growth}, 600, 800.

\bibitem[{Zhang et~al.(2000)Zhang, Chen, and Neumann}]{ZHANG2000}
Zhang, X., Chen, L., and Neumann, A.U. (2000).
\newblock The stage-structured predator–prey model and optimal harvesting
  policy.
\newblock \emph{Mathematical Biosciences}, 168(2), 201--210.

\end{thebibliography}
\end{document}